\documentclass[11pt,,xcolor=dvipsnames]{article}
\usepackage[active]{srcltx}
\usepackage{amssymb,amsmath}
\usepackage{graphicx}
\usepackage{verbatim}

\newtheorem{example}{ Example}[section]
\newtheorem{proposition}{Proposition}[section]
\newtheorem{theorem}{Theorem}[section]

\newtheorem{remark}{Remark}[section]

\numberwithin{equation}{section}

\begin{document}

\begin{center}{\large\sc
What do  `convexities' imply on Hadamard manifolds?
}\\
\vspace{0.5cm} Alexandru Krist\'aly$^*$, Chong Li$^{**}$, Genaro
Lopez$^{***}$,
Adriana Nicolae$^{****}$\\
\vspace{0.5cm}
 {$^{*}$\footnotesize  Department of Economics, Babe\c s-Bolyai
University, Cluj-Napoca, Romania\\Email address:
alexandrukristaly@yahoo.com}\\
$^{**}$ {\footnotesize  Department of Mathematics, Zhejiang
University, Hangzhou 310027, P. R. China\\Email address:
cli@zju.edu.cn}\\
$^{***}$ {\footnotesize  Departamento de An\'alisis Matem\'atico,
Universidad de Sevilla,
Apdo. 1160, 41080-Sevilla, Spain\\
             Email address: glopez@us.es}\\
$^{****}$ {\footnotesize Department of Mathematics, Babe\c s-Bolyai University, Kog\u alniceanu 1, 400084 Cluj-Napoca, Romania\\
          and\\
          Simion Stoilow Institute of Mathematics of the Romanian Academy, Research group of the project PD-3-0152, P.O. Box 1-764, RO-014700 Bucharest, Romania\\
          Email address: anicolae@math.ubbcluj.ro}
\end{center}

\vspace{1cm}

\begin{abstract}
{\footnotesize \noindent Various results based on some convexity
assumptions (involving the exponential map along with  affine maps,
geodesics and convex hulls) have been recently established on
Hadamard manifolds. In this paper we prove that these conditions are
mutually equivalent and they hold if and only if the Hadamard
manifold is isometric to the Euclidean space. In this way, we show
that some results in the literature obtained on Hadamard manifolds
are actually nothing but their well known Euclidean counterparts.}

\end{abstract}

\noindent {\it Keywords}: Hadamard manifold; convexity.\\
\noindent {\it MSC}: 53C23; 53C24.

\section{Introduction}\label{sect:1}

In recent years considerable efforts have been done to extend
concepts and results from the Euclidean/Hilbert context to settings
with no vector space structure. The motivation of such studies comes
from nonlinear phenomena which require the presence of a
non-positively curved structure for the ambient space;
 see   Jost \cite{Jost}, Krist\'aly \cite{Kristaly},
Krist\'aly, R\u adulescu and Varga \cite{KRV}, Li, L\'opez and
Mart\'\i n-M\'arquez \cite{LLMM-1}, N\'emeth \cite{Nemeth}, Udri\c
ste \cite{Udriste} and references therein.

The purpose of the present paper is to point out some conceptual
mistakes within the class of Hadamard manifolds where some authors
used equivalences between {\it convexity notions} which basically
reduce the geometric setting to the Euclidean one. Thus, in all
these papers the corresponding results and their consequences are
nothing but previously well known facts in the Euclidean case.

To be more precise, let $(M,g)$ be a Hadamard manifold (i.e., simply
connected, complete Riemannian manifold with non-positive sectional
curvature). According to the Cartan-Hadamard theorem,  the
exponential map $\exp_p:T_pM\to M$ is a global diffeomorphism for
every $p\in M$. Let $p\in M$ be fixed arbitrarily. By using the
exponential map, three convexity notions are recalled in the sequel,
mentioning also their sources without sake of completeness:
\begin{itemize}
  \item {\it Affinity.} A map $f:M\to \mathbb R$ is called affine if $f\circ
  \gamma:[0,1]\to \mathbb R$ is affine in the usual sense on $[0,1]$ for every geodesic segment  $\gamma:[0,1]\to
  M.$ Papa Quiroz \cite{Papa} and Papa Quiroz and Oliveira \cite{Papa_Oliveira} claimed that $M\ni q\mapsto g_p( \exp_p^{-1}(q),y)$ is
  affine for every $y\in T_pM$ and they used this property to prove
  convergence of various algorithms on Hadamard manifolds. This
  statement  is also used in Colao,  L\'opez,  Marino and Mart\'\i n-M\'arquez
\cite{CLMMM}, and Zhou and Huang \cite{ZH}.
  \item {\it Geodesics.} Let $q_1,q_2\in M$ be two fixed points. By construction, the unique minimal geodesic $\gamma:[0,1]\to M$ joining these points is given by
    $\gamma(t)=\exp_{q_1}(t\exp_{q_1}^{-1}(q_2)).$ Yang and Pu \cite{Yang-Pu} claimed that the curve $[0,1]\ni t\mapsto
  \exp_p((1-t)\exp_p^{-1}(q_1)+t\exp_p^{-1}(q_2))$  is
 also a minimal geodesic segment on $M$ joining the points $q_1$ and $q_2.$
  \item {\it Convex hull.}  By definition (see \cite[page 67]{Jost}), the convex hull C$(S)$ of a set $S\subset
  M$ is the smallest convex subset of $M$ containing $S.$ Instead of the convex hull, Yang and Pu \cite{Yang-Pu} introduced the {\it geodesic convex
  hull}
  GC$_p(S)$ of a set $S\subset  M$ with respect to $p\in M$ in the following way
  $${\rm GC}_p(S)=\left\{\exp_p\left(\sum_{i=1}^m\lambda_i \exp_p^{-1}(q_i)\right):\ \forall q_1,...,q_m\in S;\lambda_1,...,\lambda_m\in [0,1],\sum_{i=1}^m\lambda_i=1\right\}.$$
It is claimed  in \cite{Yang-Pu}  that  {\rm C}$(S)= {\rm GC}_p(S)$
for every $p\in M$ and $S\subset M.$
\end{itemize}

We provide below a concrete counterexample in the hyperbolic plane
which shows that the aforementioned claims are based on a
fundamental misconception.

\begin{example}\rm
Consider the Poincar\'{e} upper half-plane model $\mathbb{H} =
\{(u,v) \in \mathbb{R}^2 : v > 0\}$ endowed with the Riemannian
metric defined for every $(u,v) \in \mathbb{H}$ by
\[g_{ij}(u,v) = \frac{1}{v^2} \delta_{ij}, \quad \text{for } i,j = 1,2.\]
$(\mathbb{H},g)$ is a Hadamard manifold with constant sectional
curvature $-1$ and the geodesics in $\mathbb{H}$ are the semilines
and the semicircles orthogonal to the line $v=0$. The Riemannian
distance between two points $(u_1,v_1), (u_2,v_2) \in \mathbb{H}$ is
given by
\[d_\mathbb{H}\left((u_1,v_1), (u_2,v_2)\right) = {\rm arccosh} \left(1+\frac{(u_2-u_1)^2+(v_2-v_1)^2}{2v_1v_2}\right). \]

Fix $p=(0,1)$. By some elementary calculations (see also \cite[page
20]{Udriste}) we have that for each $(\alpha,\beta) \in T_p
\mathbb{H}$,
\[\exp_p(\alpha,\beta) = \begin{cases}
(0,e^{\beta}) &\mbox{if } \alpha = 0 \\
\left(\frac{\beta}{\alpha} + r_{\alpha,\beta}
\tanh(s_{\alpha,\beta}),\frac{r_{\alpha,\beta}}{\cosh(s_{\alpha,\beta})}
\right) & \mbox{if } \alpha \ne 0,
\end{cases} \]
where $r_{\alpha,\beta} = \sqrt{1 +
\left(\frac{\beta}{\alpha}\right)^2}$ and
\[s_{\alpha,\beta} = \begin{cases}
\sqrt{\alpha^2+\beta^2} - {\rm arcsinh} \frac{\beta}{\alpha} &\mbox{if } \alpha > 0 \\
-\sqrt{\alpha^2+\beta^2} - {\rm arcsinh} \frac{\beta}{\alpha} &
\mbox{if } \alpha < 0.
\end{cases}\]
In other words, $\exp_p(\alpha,\beta)$ belongs to the semiline
$\{u=0, v>0\}$ or a semicircle orthogonal to the line $v=0$ that
contains $p$ and for which the direction of the tangent in $p$ is
given by the vector $(\alpha,\beta)$. Moreover,
$d_{\mathbb{H}}\left(p,\exp_p(\alpha,\beta)\right) =
\|(\alpha,\beta)\| = \sqrt{\alpha^2 + \beta^2}$.

Let $q_1=\left(1,\sqrt{2}\right)$, $q_2 = \left(-1,\sqrt{2}\right)$
and take $S=\{q_1,q_2\}$. Then,
\[{\rm C}(S) = \{(u,v) \in \mathbb{R}^2 : u^2+v^2 = 3, u \in [-1,1], v > 0\}.\]

The geodesic segment joining $p$ and $q_1$ belongs to the semicircle
$K = \{(u,v) \in \mathbb{R}^2 : (u-1)^2 + v^2 = 2, v>0\}$. Denote
$\eta_1 = \exp_p^{-1}(q_1) =  (\alpha,\alpha)$ and $\eta_2 =
\exp_p^{-1}(q_2) = (-\alpha,\alpha)$ with $\alpha = \ln
(\sqrt{2}+1)^{1/\sqrt{2}}$ (see also \cite[Section 5]{WanLiYao} for
the general expression of the inverse exponential map).

Let $\eta = (1/2)\eta_1 + (1/2)\eta_2 = (0, \alpha)$ and $(0,x) =
\exp_p(\eta)$. Then $x = (\sqrt{2}+1)^{1/\sqrt{2}} > \sqrt{3}$.
Thus,  $${\rm C}(S) \ne {\rm GC}_p(S).$$ Actually, since
$(1-t)\eta_1 + t \eta_2 = \left((1-2t)\alpha,\alpha\right)$, the set
${\rm GC}_p(S)$ consists of the points
$\exp_p\left((1-2t)\alpha,\alpha\right)$, with $t \in [0,1]$, see
Figure \ref{fig:1}.

\begin{figure}[h!]
\begin{center}
\includegraphics[scale=0.48]{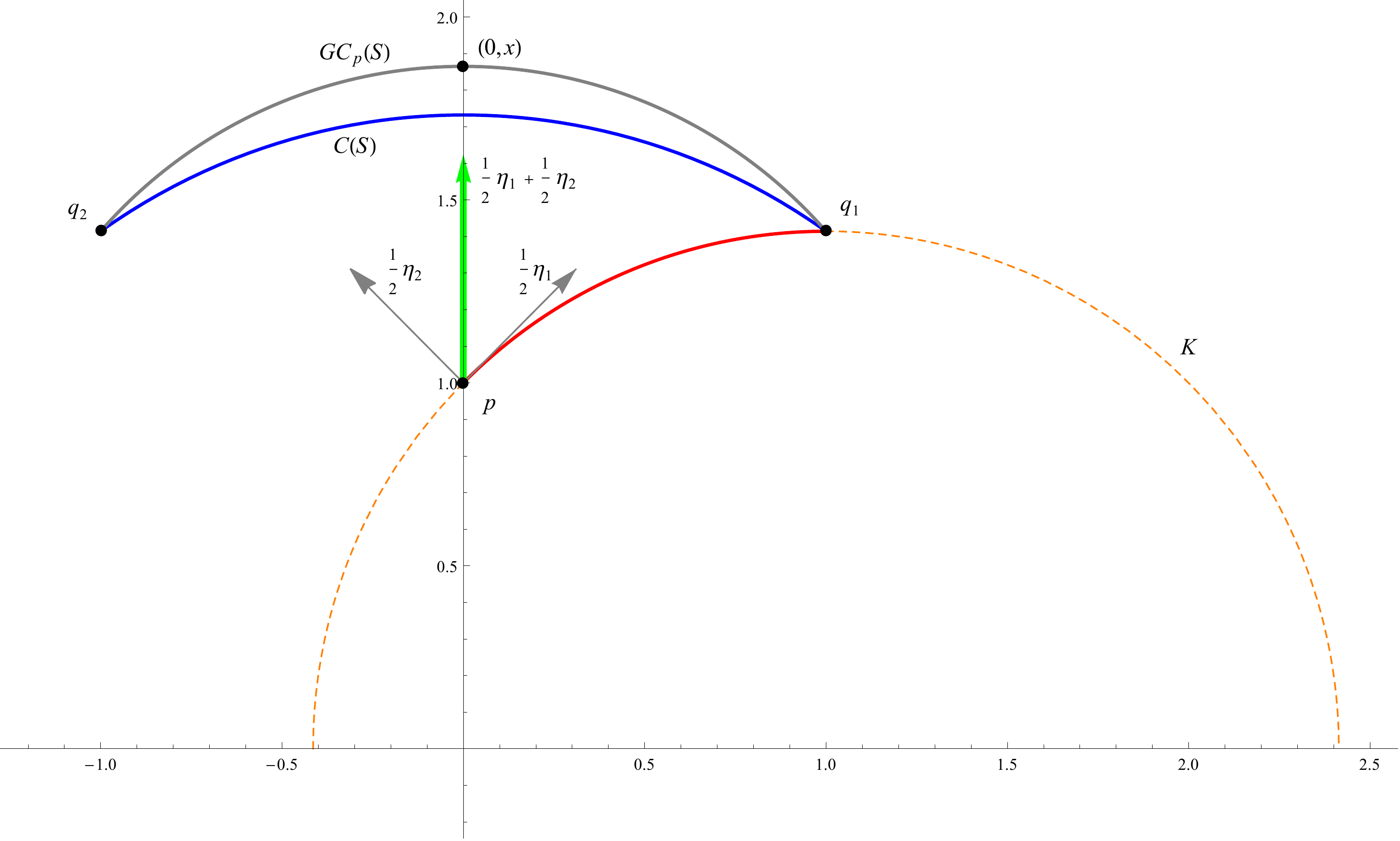}
\caption{Difference between convex hull and geodesic convex hull.}
\label{fig:1}
\end{center}
\end{figure}

Note that the above fact also shows that the curve $[0,1]\ni
t\mapsto \exp_p((1-t)\exp_p^{-1}(q_1)+t\exp_p^{-1}(q_2))$ is not the
minimal geodesic joining $q_1$ and $q_2$.

It is well known that any affine mapping defined on $\mathbb{H}$ is
constant. In particular, one can choose $y \in T_p
\mathbb{H}\setminus \{0\}$ such that the mapping $\mathbb{H}\ni
q\mapsto \langle \exp_p^{-1}(q),y \rangle$ is not constant and so it
is not affine. This fact can be checked directly as well. Moreover,
it is also obvious that $\exp_p:T_p \mathbb{H}\to \mathbb{H}$ is not
a global isometry.
\end{example}

\section{Main result}\label{sect:22}

As the main result of this paper we prove the following rigidity
theorem.


\begin{theorem}\label{theorem-1}
Let $(M,g)$ be an $n-$dimensional Hadamard manifold and $p\in M.$
Then the following statements are equivalent:
\begin{itemize}
  \item[{\rm (i)}] The map $M\ni q\mapsto g_p( \exp_p^{-1}(q),y)$ is
  affine for every $y\in T_pM;$
  \item[{\rm (ii)}] For every $q_1,q_2\in M$, the curve $[0,1]\ni t\mapsto
  \exp_p((1-t)\exp_p^{-1}(q_1)+t\exp_p^{-1}(q_2))$  is
  the minimal geodesic segment joining the points $q_1$ and $q_2;$
  \item[{\rm (iii)}] For every non-empty set $S\subset M$,
  {\rm C}$(S)={\rm GC}_p(S);$
  \item[{\rm (iv)}] The map $\exp_p:T_pM\to M$ is a global isometry$;$
\item[{\rm (v)}] The sectional curvature on $(M,g)$ is identically zero $($i.e., $(M,g)$
is isometric to the usual Euclidean space $(\mathbb R^n,e)).$
\end{itemize}
\end{theorem}

In order to prove Theorem \ref{theorem-1}, we recall two results.

\begin{proposition}{\rm [Choquet theorem; see \cite[Theorem 6.5]{Udriste}]}\label{Choquet}
 An $n-$dimensional
Riemannian manifold $(M,g)$ is the Riemannian product of an
$(n-p+1)-$dimensional Riemannian manifold and the Euclidean space
$\mathbb R^{p-1}$ $($at least locally$)$ if and only if the vector
space of all affine functions on $M$ has dimension $p.$ $[$In
particular, the sectional curvature restricted to the components of
$\mathbb R^{p-1}$ is identically zero.$]$
\end{proposition}

\begin{proposition}{\rm [See \cite[Lemma 3.3.1]{Jost}]}\label{convex-hull-character}
The convex hull {\rm C}$(S)$ of a set $S\subset M$ is
$${\rm C}(S)=\bigcup_{k=0}^\infty S_k,$$
where $S_0=S$ and for every $k\in \mathbb N$, $S_k$ is the union of
all geodesic segments between points of $S_{k-1}.$
\end{proposition}

\noindent Now, we are ready to prove our rigidity result.\\

{\it Proof of Theorem \ref{theorem-1}.} The equivalence
(iv)$\Leftrightarrow$(v) is trivial, see \cite[Theorem
4.1]{doCarmo}.

(v)$\Rightarrow$(ii) This implication follows directly because
property (ii) is satisfied in the Euclidean space $(\mathbb R^n,e)$
and geodesics are invariant by isometries between Hadamard
manifolds.

(ii)$\Rightarrow$(iii) is also trivial, coming from the two
definitions and elementary computations.

(ii)$\Rightarrow$(i) Let $y\in T_pM$ be fixed arbitrarily; for
convenience, let $f_y:M\to \mathbb R$ be defined by $f_y(q)=g_p(
\exp_p^{-1}(q),y)$. By  assumption, any geodesic segment in $M$ can
be represented by $\gamma=\exp_p\circ \gamma_0$, where
$\gamma_0(t)=(1-t)\exp_p^{-1}(q_1)+t\exp_p^{-1}(q_2)$, $t\in [0,1]$
for some $q_1,q_2\in M.$ Then,
$$f_y(\gamma(t))=g_p( \exp_p^{-1}(\gamma(t)),y)=g_p(
\gamma_0(t),y)=(1-t)g_p(\exp_p^{-1}(q_1),y)+tg_p(\exp_p^{-1}(q_2),y),$$
which is an affine function on $[0,1]$ in the usual sense.

(i)$\Rightarrow$(v) We show that the dimension of  the space of
affine functions on $M$ is $n+1.$ By assumption, $f_y(q)=g_p(
\exp_p^{-1}(q),y)$ is an affine function on $M$ for every $y\in
T_pM$.  In particular, it follows that ${\rm Hess}_gf_y=0,$ since
$f_y$ is both convex and concave, see \cite{Udriste}. Since
Hess$_gf_y(V,W)=g(\nabla_V {\rm grad}f_y,W)$ for every vector fields
$V,W$ on $M$, the latter relation implies in particular that
grad$f_y$ is a parallel vector field along any geodesic of $M$.
Since dim$(M)=n$, we may fix $y_1,...,y_n\in T_pM$ such that in
every $q\in M$, the set $\{{\rm grad}f_{y_1},...,{\rm
grad}f_{y_n}\}$ forms a basis of the tangent space $T_qM$
(basically, it is enough to guarantee this property just in one
point and use parallel transport at any fixed point). In this
manner, we constructed exactly $n$ non-constant, linearly
independent affine functions $f_{y_i}$ on $M$, $i=1,...n,$
corresponding to the elements $y_1,...,y_n\in T_pM$. Moreover, we
may add to this set also a constant function which is affine and
linearly independent of $\{f_{y_1},...,f_{y_n}\}.$ Therefore, the
vector space of all affine functions on $M$ has dimensional  $p\geq
n+1$. According to Choquet theorem (see Proposition \ref{Choquet}),
we also have that $p\leq n+1$. Therefore, $p=n+1$ and again by
Choquet theorem and from the fact that $(M,g)$ is a Hadamard
manifold, it follows that $M$ is globally represented as $\mathbb
R^n$ whose sectional curvature
 is identically zero.

(iii)$\Rightarrow$(ii) Let $q_1,q_2\in M$ and  $S=\{q_1,q_2\}$. On
the one hand, by the definition of the geodesic convex hull, since
$S$ contains just two elements, we clearly have that
\begin{eqnarray*}
  \text{GC}_p(S) &=& \left\{\exp_p\left(\lambda_1 \exp_p^{-1}(q_1)+\lambda_2
\exp_p^{-1}(q_2)\right):\lambda_1,\lambda_2\in [0,1],\
\lambda_1+\lambda_2=1\right\} \\
   &=& \left\{\exp_p\left((1-t) \exp_p^{-1}(q_1)+t
\exp_p^{-1}(q_2)\right):t\in [0,1]\right\}.
\end{eqnarray*}
On the other hand, if $S_0=S$, then the set $S_1$ in Proposition
\ref{convex-hull-character} is precisely the image of the unique
minimal geodesic segment $\gamma:[0,1]\to M$ joining the points
$q_1$ and $q_2$. Moreover, if we take any two points in $S_1={\rm
Im}(\gamma)$ and join them by a geodesic segment, the minimality of
$\gamma$ implies that the image of the latter geodesic will be a
subset of ${\rm Im}(\gamma)$. Therefore, $S_1=S_2=...={\rm
Im}(\gamma)$. Consequently, ${\rm C}(S)={\rm Im}(\gamma).$ Since by
assumption ${\rm C}(S)=\text{GC}_p(S),$ one obtains (ii).
 \hfill $\square$

\begin{remark}\rm
Further definitions and open questions concerning the convex hull on
non-positively curved spaces can be found in the literature, see
Ledyaev, Treiman and Zhu \cite[Conjecture 1]{LTZ} and Nava-Yazdani
and Polthier \cite{Nav-Pol}. It would be interesting to study the
relationship between these notions and the geometric structure of
the ambient space.
\end{remark}

\vspace{0.5cm} \noindent {\bf Acknowledgment.} A. Krist\'aly was
supported by a grant of the Romanian Ministry of Education,
CNCS-UEFISCDI, project number PN-II-RU-TE-2011-3-0047. C. Li and G.
L\'opez were supported by DGES (Grant MTM2012-34847-C02-01). A.
Nicolae was supported by a grant of the Romanian Ministry of
Education, CNCS-UEFISCDI, project number PN-II-RU-PD-2012-3-0152.
Part of this work was carried out while some of the authors were
visiting the University of Seville supported by DGES (Grant
MTM2012-34847-C02-01).

\end{document}